# A NOTE ON NOETHERIAN RINGS.

C.L.Wangneo

Jammu,J&K,India,180002

(E-mail:-wangneo.chaman@gmail.com )

**Abstract:-** In this paper we introduce the definition of a noetherian disjoint ring and that of a noetherian non-disjoint ring . For a noetherian ring R , with nilradical N if P and Q represent the semiprime ideals of R called as the right and the left krull-homogenous parts of N as defined in [8] , then we prove the main theorem of this paper for the ring R whose statement is given below.

**Main Theorem :-** Let R be a Noetherian ring with nilradical N . Let P and Q represent the right and the left krull-homogenous parts of N . Then the following hold true for the ring R ;

(a) If R is a disjoint ring , then the nilradical N of R is a right and a left weakly ideal invariant ideal of R . Hence N is a right and a left localizable semiprime ideal of R .

(b) If R is a non-disjoint ring then the following are equivalent conditions on R ;

 (i) N is a right and a left weakly ideal invariant ideal of R .

 (ii) P = Q is a right and a left localizable semiprime ideal of R .

**Introduction :-** In this paper we introduce the definition of a noetherian disjoint ring and that of a noetherian non-disjoint ring . For a noetherian ring R , with nilradical N if P and Q represent the semiprime ideals of R called as the right and the left krull-homogenous parts of N as defined in [8] , then we prove the main theorem of this paper, namely theorem (3.8) for the ring R which states the following.

**Theorem (3.8) :-** Let R be a Noetherian ring with nilradical N. Let P and Q represent the right and the left krull-homogenous parts of N. Then the following hold true for the ring R;

(a) If R is a disjoint ring, then the nilradical N of R is a right and a left weakly ideal invariant ideal of R. Hence N is a right and a left localizable semiprime ideal of R.

(b) If R is a non-disjoint ring then the following are equivalent conditions on R;

(i) N is a right and a left weakly ideal invariant ideal of R.

(ii) P = Q is a right and a left localizable semiprime ideal of R.

The paper is divided into three sections. In section (1) we introduce the preliminaries for this paper. We thus recall the following definitions of ideals A(R) and B(R) of a noetherian ring R, where A(R) is the sum of all right Ideal I of R such that $|I|_r < |R|_r$, and B(R) is the sum of all left ideal K of R with $|K|_l < |R|_l$. We then state some some basic results regarding the ideals A(R) and B(R). Next we recall for a noetherian ring R certain subsets $\Lambda(R)$ and $\Lambda'(R)$ of min. spec.R as follows;

$\Lambda(R) = \{ P_i \text{ in spec. }(R) / |R/P_i|_r = |R|_r \}$ and

$\Lambda'(R) = \{ Q_i \text{ in spec. }(R) / |R/Q_i|_l = |R|_l \}$.

These definitions allow us to consider the semiprime ideals P and Q of R, where $P = \cap p_i$; $p_i \in \Lambda$ and $Q = \cap q_i$; $q_i \in \Lambda'$. We call P and Q as the right and the left krull-homogenous parts of N as defined in [8]. We come up with two obvious cases for the noetherian ring R, namely,

**Case (i) :-** When $\Lambda \cap \Lambda' = \Phi$, and

**Case (ii) :-** When $\Lambda \cap \Lambda' \neq \Phi$.

We then characterize in section(1) the above two cases by conditions equivalent to them. Finally we end section (1) with the following definitions. We call a noetherian ring R a disjoint ring if $\Lambda \cap \Lambda' = \Phi$ and we call a noetherian ring R a non-disjoint ring if $\Lambda \cap \Lambda' \neq \Phi$.

In section (2) we study a noetherian disjoint ring and in section(3) we study a noetherian non-disjoint ring. We prove our key theorems of these sections, namely, theorems (2.9) and (3.7), respectively. Finally in section (3) we combine the key theorems of sections (2) and (3) as theorem (3.8) which we call as our main theorem.

<u>**Notation and Terminology**</u>:-

Throughout in this paper a ring is meant to be an associative ring with identity which is not necessarily a commutative ring. We will throughout adhere to the same notation and terminology as in [3] and [8]. Thus by a noetherian ring R we mean that R is both a left as well as a right noetherian ring. By a module M over a ring R we mean that M is a right R-module unless stated otherwise. For the basic definitions regarding noetherian modules and noetherian rings and all those regarding Krull dimension, we refer the reader to [3]. Moreover we will use the following terminology throughout. If R is a ring then we denote by Spec.R, the set of prime ideals of R and by min.Spec.R, the set of minimal prime ideals of R. Again if R is a ring and M is a right R module, then we denote by r(T) the right annihilator of a subset T of M and by l(T) the left annihilator of a subset T of W in case W is a left R module. Also $|M|_r$ denotes the right Krull dimension of the right R-module M if it exists and $|W|_l$ denotes the left Krull dimension of the left R module W if it exists. For two subsets A and B of a given we denote by A ≤ B that B contains A and A < B denotes A ≤ B but A≠B. Also for two sets A and B, A⊄B denotes that the subset B does not contain the subset A. For an ideal A of R, c(A) denotes the set of elements of R that are regular modulo A.

**Section (1) (Preliminaries)** :- In this section we first recall the following definitions from [ 8 ].

**Definition and Notation(1.1)** :- Let R be a noetherian ring. Denote by A(R) the sum of all right Ideal I of R such that $|I|_r < |R|_r$. Similarly denote by B(R) the sum of all left ideal K of R with $|K|_l < |R|_l$. If there is no confusion regarding the underlying ring R we may write A and B for A(R) and B(R) respectively.

We now state the following basic results regarding the above definitions from [ 8 ] ;

**Proposition (1.2)** :- Let R be a Noetherian ring. Let A and B be as in definition(1.1) above. Then the following hold true ;

(i) A is an ideal of R and is the unique largest right ideal of R with

  $|A|_r < |R|_r$ and $|R/A|_r = |R|_r$.

(ii) Moreover R/A is a right k-homogenous ring and if $A \neq 0$, then

  $r(A) \neq 0$ and $|R/r(A)|_r < |R|_r$.

**Proof** :- For the proof of the above result we ask the reader to see the proof of the proposition (1.3 ) of [ 8 ].

We now state proposition (1.3) below which is similar to proposition (1.2) above and may be considered as its left analogue ;

**Proposition(1.3)** :- Let R be a Noetherian ring. Let A and B be as in definition(1.1) above. Then the following hold true ;

(i) B is an ideal of R and is the unique largest left ideal of R with $|B|_l < |R|_l$ and $|R/B|_l = |R|_l$.

(ii) Moreover R/B is a left k-homogenous ring and if $B \neq 0$, then $l(B) \neq 0$ and $|R/l(B)|_l < |R|_l$.

Next we recall the following definition and notation from [8].

**Definition and Notation (1.4)** :- (i) Let R be a Noetherian ring with nilradical R. Define subsets of min. spec. R denoted by $\Lambda(R)$ and $\Lambda'(R)$ as follows ;

$\Lambda(R) = \{ P_i \text{ in spec. } (R) / |R/P_i|_r = |R|_r \}$ and

$\Lambda'(R) = \{ Q_i \text{ in spec. } (R) / |R/Q_i|_l = |R|_l \}$.

We define another subset namely X(R) of min.spec.(R) as follows ;

$X(R) = \{ l \in \text{min.Spec.}(R) / |R/l|_r < |R|_r, \text{ and } |R/l|_l < |R|_l \}$

If there is no confusion about the underlying ring R we set $\Lambda = \Lambda(R)$, $\Lambda' = \Lambda'(R)$ and $X = X(R)$.

(ii) It is not difficult to see that if $P_i \in \Lambda(R)$, then $P_i \in$ min. spec.(R). We denote by P(R) the ideal $P(R) = \cap P_i$, $P_i \in \Lambda$ and by Q(R) the ideal $Q(R) = \cap Q_i$, $Q_i \in \Lambda'$. Then the ideals P(R) and Q(R) are called the right and the left krull homogenous parts respectively of the nilradical N(R) of the ring R. If there is no confusion regarding the underlying ring R we may write P for P(R) and Q for Q(R). We may note that in this case the factor rings R/P and R/Q are semiprime right and left krull homogenous rings respectively. We call X the non krull-homogenous part of N. It is also clear that $N \leq P$ and $N \leq Q$.

**(iii) For an ideal I of R we have ;  $\Lambda(R/I) = \{ P_i/I$ in spec. $(R/I) / I \leq P_i$ and $|R/P_i|_r = |R/I|_r \}$  and  $\Lambda'(R/I) = \{ Q_i/I$ in spec. $(R/I) / I \leq Q_i$  and $|R/Q_i|_l = |R/I|_l \}$ .**

**We now recall definition (2.3) of [8] of a right ( or a left ) weakly krull symmetric ring , a weakly krull symmetric  and a weakly ideal krull symmetric noetherian ring  R and call it as  our definition (1.5) given below .**

**<u>Definition (1.5)</u> :- Let R be a Noetherian  ring  with nilradical N and let  $\Lambda$ , $\Lambda'$ , A , B be as  defined in definitions (1.4 ) and (1.1) above . Let   P =  $\cap$ $p_i$ ;  $p_i$ $\mathcal{E}$ $\Lambda$  and let  Q =  $\cap$ $q_i$ ;  $q_i$ $\mathcal{E}$ $\Lambda'$ . Then  we  have the following definitions ,**

**(i) R is called a  right ( or  a  left) weakly krull symmetric ring if $\Lambda'(R) \leq \Lambda(R)$ ( or  $\Lambda(R) \leq \Lambda'(R)$) .**

**(ii) R is called  a  weakly krull symmetric ring if  $\Lambda(R) = \Lambda'(R)$ .**

**(iii) R is called  a  weakly ideal  krull symmetric ring if  A = B .**

**<u>Definition (1.6) ( Two cases )</u> :- We now   state and  define two obvious cases for a noetherian ring R , namely ,**

**<u>Case  (i) ( the  disjoint case )</u> :- When  $\Lambda \cap \Lambda' = \Phi$, and**

**<u>Case (ii) ( the  non-disjoint  case )</u> :-  When  $\Lambda \cap \Lambda' \neq \Phi$ .**

**We characterise the  above  two  cases  by   conditions  equivalent  to them in the theorems  stated  below ;**

**<u>Theorem (1.7)</u> :- Let R be a  Noetherian  ring  with  nilradical N and let $\Lambda$ , $\Lambda'$ , A , B be as  defined in definitions (1.4) and (1.1) above . Then the following are equivalent ;**

**(a)   $\Lambda \cap \Lambda' = \Phi$ .**

**(b) $|R/A|_l < |R|_l$ .**

(c) $|R/B|_r < |R|_r$.

**Proof :-** We first prove (a) => (b). For this asume (a), then we show that (b) is true. Suppose not and let $|R/A|_l = |R|_l$. Then there exists a prime ideal p, such that $A \leq p$ and $|R/p|_l = |R|_l$. Thus $p \, \varepsilon \, \Lambda'(R)$ and so $p \, \varepsilon \, \Lambda'(R/A)$. Since R/A is right krull homogenous, so by the main theorem of [8], namely theorem (2.10) of [8], we have that $\Lambda'(R/A) \leq \Lambda(R/A)$. Thus we must have that $p/A \, \varepsilon \, \Lambda(R/A)$. Hence $|R/p|_r = |R|_r$. Thus $p \, \varepsilon \, \Lambda(R)$ also. Hence $\Lambda \cap \Lambda' \neq \Phi$, which is a contradiction to (a). Hence we must have that $|R/A|_l < |R|_l$.

(b) => (a) We now prove (b) => (a). If $|R/A|_l < |R|_l$, then we show that $\Lambda \cap \Lambda' = \Phi$. Suppose not, then $\Lambda \cap \Lambda' \neq \Phi$. Let p be a prime ideal such that $p \, \varepsilon \, \Lambda \cap \Lambda'$, then using propositions (1.2) and (1.3) above we get that $A \leq p$, and $B \leq p$. Since $|R/p|_l = |R|_l$, this must imply that $|R/A|_l = |R|_l$ which contradicts (b). Hence we must have that $\Lambda \cap \Lambda' = \Phi$. This proves (a) <=> (b).

(a) <=> (c) The proof of (a) <=> (c) is similar to the proof of (a) <=> (b) given above.

A reformulation of the above theorem (1.7) is theorem (1.8) stated below ;

**Theorem (1.8)** :- Let R be a Noetherian ring with nilradical N and let $\Lambda, \Lambda', A, B$ be as defined in definitions (1.4) and (1.1) above. Then the following are equivalent ;

(a) $\Lambda \cap \Lambda' \neq \Phi$.

(b) $|R/A|_l = |R|_l$ (it is always true that $|R/A|_r = |R|_r$).

(c) $|R/B|_r = |R|_r$ (it is always true that $|R/B|_l = |R|_l$).

**Deinition and Notation (1.9)** :- If $\Lambda \cap \Lambda' \neq \Phi$ , then set V(R) = $\Lambda$(R) $\cap$ $\Lambda'$ (R) and denote by Y(R) the ideal Y(R) = $\cap$ x / { x $\varepsilon$ $\Lambda$(R) $\cap$ $\Lambda'$ (R)} . Then we have the following theorem regarding the case (2) mentioned above .

**Theorem (1.10)** :- Let R be a Noetherian ring with nilradical N and let $\Lambda$ , $\Lambda'$ , A , B be as defined in definitions (1.4) and (1.1) above . If $\Lambda \cap \Lambda' \neq \Phi$ , then the following hold true ;

(a) Y is a semiprime ideal of R such that A $\leq$ Y and B $\leq$ Y . Hence A+B $\leq$ Y .

(b) Moreover V ( R/A+B) = $\Lambda$ (R/A+B) = $\Lambda'$ (R/A+B ) .

**Proof** :- (a) Using propositions (1.2) and (1.3) respectively it is not difficult to see that Y is a semiprime ideal of R such that A $\leq$ Y and B $\leq$ Y . Hence A+B $\leq$ Y .

(b) (b) follows easily from the definition of the ideals A , B and subsets $\Lambda$ ( R/A+B) and $\Lambda'$ ( R/A+B) of min.spec.(R/A+B) .

We finally end section (1) with the following definition .

**Definition (1.11)** :- Let R be a noetherian ring . If A , B and $\Lambda$ , $\Lambda'$ are as defined in definitions (1.1) and (1.4) respectively , then we have the following definitions ;

(a) We call R a disjoint ring if $\Lambda \cap \Lambda' = \Phi$ .

(b) We call R a non-disjoint ring if $\Lambda \cap \Lambda' \neq \Phi$ .

**Section (2) (Noetherian disjoint rings )** :- In this section we study a noetherian disjoint ring and we prove our key theorem of this section , namely theorem (2.9) . To do so we must recall the definitions of the notion of the weak ideal invariance (w.i.i for short ) and that of the localizability of an ideal of a noetherian ring R as stated for example in [2] . This we do below ;

**Definition (2.1) (Weak ideal invariance) :-** Let R be a noetherian ring. If I is any ideal of R we call I a right weakly ideal invariant ideal ( right w.i.i for short ) if for any right ideal J of R with $|R/J|_r < |R/I|_r$, we have that $|I/JI|r < |R/I|r$. Similarly we define the notion of the left w.i.i of an ideal I of a noetherian ring R. An ideal is said to be w.i.i if it is left w.i.i as well as right w.i.i.

**Definition (2.2)** ( localisability ) :- Let I be an ideal of a ring R. Let $C(I) = \{ c \varepsilon R / c+I$ is regular in R/I $\}$. The ideal I is said to be right localisable if the elements of C(I) satisfy the right ore condition ; that is, given $c \varepsilon C(I)$ and $x \varepsilon R$, there exists $d \varepsilon C(I)$ and $y \varepsilon R$ such that $xd = cy$. An ideal is said to be localisable if it is left and right localisable. Define $T(I) = \{ x \varepsilon R / xc=0$, for some $c \varepsilon C(I) \}$. Note that, if I is right localisable, then T(I) is actually an ideal.

**Remark :-** We now state theorem (2.3) below which states that the nilradical N of a noetherian ring R is a right and a left localizable semiprime ideal of R if we assume that N is a right and a left w.i.i ideal of R.

**Theorem (2.3) :-** Let R be a Noetherian ring with nilradical N. If N is right and left weakly ideal invariant ideal of R, then N is a right and a left localizable semiprime ideal of R.

**Proof :-** We assume that N is a right and a left w.i.i semiprime ideal of R. We then show that N is a right and a left localizable semiprime ideal of R. To see this we make the following claim first ;

**Claim (1)** :- For any $d \in c(N)$, $|R/dR+N|_r < |R|_r$ if and only if $|R/dR|_r < |R|_r$.

**Proof** :- To prove the above claim we first observe that N is a nilpotent ideal of R. Now given that N is a right w.i.i ideal of R, so apply lemma (3) of [5] to conclude that $|R/dR+N|_r < |R|_r$ implies that $|R/dR|_r < |R|_r$. The converse is clearly easy.

Now use the fact that N is a right w.i.i ideal of R together with the proof outlined in the beginning of theorem (8) of [5] to conclude that N is a right localizable semiprime ideal of R. That N is a left localizable semiprime ideal follows similarly.

**Remark:-** To the best of the authors knowledge the converse of the above theorem is an open question.

We now recall briefly from [3] the following definitions ;

**Definition (2.4)** :- (a) For a finitely generated right R-module M over a noetherian ring R we denote by Ass.M the set of prime ideals of R associated to M on the right. We use a similar notation for a left R-module M causing no confusion.

(b) Again from [3] we call a right (or a left) R-module M to be right (or left) primary if it has a unique right (or left) associated prime ideal.

We now recall further the following definition from [1] ;

**Definition (2.5)** :- A noetherian ring R is said to have the right large condition if for any essential right ideal J of R we have that $|R/J|_r < |R|_r$

With these definitions we now state the following proposition ;

**Proposition (2.6)** :- Let R be a Noetherian ring. Let p be a prime ideal of R such that $|R/p|_r = |R|_r$. Let M be a finitely generated critical right R-module M with $|M|_r = |R|_r$ and Ass.(M) = p. If r(M) = I, then we have the following ; (i) $I \leq p$.

(ii) R/I is a right krull homogenous ring with $|R/I|_r = |R|_r$.

(iii) R/I is a p/I primary ring such that p/I is a non-essential right ideal of R/I and moreover R/I satisfies the right large condition.

(iv) Either I=p or if $I \neq p$, then there exists a nonzero ideal $J \neq p$, I < J, with $Jp \leq I$ and $|R/J|_r < |R|_r$.

**Proof** :- Use [4], lemma (1.9) and [1] corollary (2.9) for the proof of (i), (ii) and (iii) respectively. We will now prove (iv). To prove (iv) observe that since R/I is a p/I primary ring and since p/I is non-essential as a right ideal of R/I thus there must exist an ideal J of R with I < J and $Jp \leq I$. Now two cases arise ;

case (i) Either J=R in which case I=p.

case (ii) In this case J/I is a proper large right ideal of R/I and hence we must have that $|R/J|_r < |R|_r$.

**Remark:-** The above proposition immediately yields the important lemma (2.7) below which is what we will need throughout this section ;

**Lemma (2.7)** :- Let R be a Noetherian ring. Let p be a prime ideal of R such that $|R/p|_r = |R|_r$ and $|R/p|_l < |R|_l$. Then for any finitely generated critical right R-module M with $|M|_r = |R|_r$ and Ass.(M) = p, we have that r(M) = p.

**Proof** :- Let r(M) = I. Clearly $I \leq p$. From proposition (2.6) above, we have the following ,

If $I \neq p$, then there exists a nonzero ideal $J \neq p$, $I < J$, and $Jp \leq I$, and $|R/J|_r < |R|_r$. Now since $R/I$ is a right krull homogenous ring, thus from theorem (2.10) of [8] we get that $\Lambda'(R/I) \leq \Lambda(R/I)$. Hence we must have that $|R/J|_l < |R|_l$. Since $Jp \leq I$, this yields that $|p|_l < |R|_l$, and thus $|R/p|_l = |R|_l$, a contradiction to the given hypothesis. Hence, we must have that $I = p$.

We now state Theorem (2.8) below which gives conditions equivalent to the w.i.i of the nilradical of a noetherian ring $R$.

<u>Theorem(2.8)</u> :- Let $R$ be a noetherian ring with nilradical $N$. Then the following conditions are equivalent ;

(i) $N$ is right w.i.i ideal.

(ii) For any prime ideal $q$ of $R$ with $|R/q|_r = |R|_r$, $q$ is right w.i.i.

(iii) For any ideal $I$ of $R$ with $|R/I|_r = |R|_r$, $I$ is right w.i.i.

(iv) If $M$ is any finitely generated critical right R-module with $|M|_r = |R|_r$ and $Ass.(M) = q$, then $r(M) = q$.

<u>Proof</u> :- This is proved as theorem (2.5) of [2].

We now prove the key theorem of this section below.

<u>Theorem (2.9)</u> :- Let $R$ be a Noetherian ring with nilradical $N$. If $R$ is a disjoint ring then $N$ is a right and a left weakly ideal invariant ideal of $R$. Hence $N$ is a right and a left localizable semiprime ideal of $R$.

<u>Proof</u> :- Let $\Lambda$, $\Lambda'$ be as defined in definition (1.4) of section (1) above. If $R$ is a disjoint ring, then we have that $\Lambda \cap \Lambda' = \Phi$. We will first show that $N$ is a right w.i.i ideal of $R$. Let $M$ be a finitely generated critical right R-module with $|M|_r = |R|_r$. If $Ass.(M) = p$, then clearly $|R/p|_r = |R|_r$. Since $\Lambda \cap \Lambda' = \Phi$, thus we must have that $|R/p|_l < |R|_l$. Thus from the above lemma (2.7), we get that $r(M)$

= p . Hence using theorem (2.8) stated above, we get that N is right w.i.i ideal of R. Similarly we can show that N is a left w.i.i ideal of R. That N is a right and a left localizable semiprime ideal of R follows from theorem (2.3) stated above.

**Corollary (2.10) :-** Let R be a Noetherian ring with nilradical N (R). Let R/I be a factor ring of R with $|R/I|_t = |R|_r$. Then we have either

(i) $|R/I|_l = |R|_l$ or

(ii) $|R/I|_l < |R|_l$, in which case we must have that N (R/I) is right weakly ideal invariant ideal of R/I .

**Proof :-** If (i) is true, then there is nothing to prove.

So assume (ii), namely, that $|R/I|_l < |R|_l$. We show in this case that N(R/I) is a right w.i.i ideal of R/I. To see this let M be a finitely generated, critical right R/I- module with $|M|_r = |R/I|_r = |R|_r$. Let Ass.(M) = p/I, for some prime ideal p of R such that I ≤ p. Clearly $|R/p|_r = |R|_r$ and $|R/p|_l < |R|_l$. Now obviously, we can consider M as a finitely generated, critical, right R-module, with Ass.(M) = p and $|M|_r = |R|_r$. From lemma (2.7) stated above we get that r(M) = p . Now this is true for any finitely generated, critical right R/I module with $|M|_r = |R/I|_r$. Hence from theorem (2.8) above, we must have that N(R/I) is a right w.i.i ideal of R/I .

**Section (3) (Noetherian non-disjoint rings ) :-** In this section we study a noetherian non-disjoint ring. We will prove our key theorem of this section, namely, theorem (3.7), whose statement is given below ;

**Theorem (3.7) :-** Let R be a noetherian ring with nilradical N . If R is a non-disjoint ring then the following are equivalent conditions on R ;

(a) N is right and left w.i.i. ideal of R.

(b) P = Q is both a right and left localizable ideal of R.

To prove this theorem we recall the following notation from section (1) above.

**Notation (3.1)** :- (i) For a noetherian ring R with nilradical N recall from section (1) above that $\Lambda$ and $\Lambda'$ are the subsets of min.spec.R as defined in definition (1.4).

(ii) Next for the noetherian ring R we recall again from section (1) above the definitions of the ideals P and Q of R which we call as the right and the left krull homogenous parts respectively of the nilradical N of the ring R.

**Remark:-** We now state theorem (3.2) below which is the main theorem of [9], namely, theorem (2.10) of [9].

**Theorem (3.2)** :- Let R be a Noetherian ring with nilradical N and let $\Lambda$, $\Lambda'$, A, B be as defined in definitions (1.4) and (1.1) of section (1) above. Let $P = \cap p_i$; $p_i \in \Lambda$ and $Q = \cap q_i$; $q_i \in \Lambda'$. Then the following hold true ;

(a) N is a right w.i.i ideal of R if and only if P is a right localizable ideal of R.

(b) N is a left w.i.i ideal of R if and only if Q is a left localizable ideal of R.

**Remark** :- To continue with the proof of our main theorem we recall that we have already stated two obvious cases above for a noetherian ring R, namely,

Case (i) :- When $\Lambda \cap \Lambda' = \Phi$, and

Case (ii) :- When $\Lambda \cap \Lambda' \neq \Phi$.

**We now discuss case (2) ;**

**Definition and notation (3.3) :-** Let R be a noetherian ring with nilradical N. Let $\Lambda, \Lambda'$ be as in definition (1.4) of section (1) above. If R is a non-disjoint ring, that is if $\Lambda \cap \Lambda' \neq \Phi$, then set $V = \Lambda \cap \Lambda'$ and let $Y = \{ \cap p \ / \ p \ \varepsilon \ V \}$. With these notation we have the following theorem ;

**Theorem(3.4) :-** Let R be a Noetherian ring with nilradical N and let $\Lambda, \Lambda', A, B$ be as defined in definitions (1.4) and (1.1) respectively. Let $P = \cap p_i$ ; $p_i \ \varepsilon \ \Lambda$ and let $Q = \cap q_i$ ; $q_i \ \varepsilon \ \Lambda'$. Then the following are true ;

(i) If R is right krull homogenous, then $\Lambda' \leq \Lambda$, and $P \leq Q$.

(ii) If R is right krull homogenous and N is right w.i.i ideal of R, then $N = P$.

**Proof :-** (i) Since R is right krull homogenous, so apply the main theorem of [8] to conclude that $\Lambda' \leq \Lambda$. Clearly then $P \leq Q$.

(ii) If R is right krull homogenous and N is right w.i.i ideal of R then we prove that $N = P$.

To prove this we first make the following claim ;

**Claim(1) :-** Let $q, p$ be minimal prime ideals of R, with $|R/p|_r = |R|r$, and $|R/q|_r < |R|r$. Then there exists an ideal B with $|R/B|_r < |R|r$, such that $pB \leq qp$.

**Proof :-** This is true because N is right w.i.i ideal of R.

**Proof of (ii) :-** We now prove (ii). To prove (ii) assume

$A_1 A_2 A_3 \ldots A_m \ldots A_n = 0$, where $A_i$ are all the minimal prime ideals of R. Then applying claim (1) we get that $p_1 p_2 \ldots p_k D = 0$, where each $p_i \ \varepsilon \ \Lambda$ and D is an ideal such that $|R/D|_r < |R|r$. Hence since R is right krull homogenous we must have that P is a nilpotent ideal of R. Thus N = P.

**Theorem(3.5)** :- Let R be a Noetherian ring with nilradical N and let $\Lambda$, $\Lambda'$, A, B be as defined in definitions (1.4) and (1.1) respectively. Let $P = \cap p_i$; $p_i \in \Lambda$ and let $Q = \cap q_i$; $q_i \in \Lambda'$. If R is right krull homogenous and N is right w.i.i ideal of R, then R is also a left krull homogenous ring.

**Proof** :- To prove the theorem we first observe from theorem (3.4)(i) above that since R is a right krull homogenous ring, hence $\Lambda' \leq \Lambda$ and thus $P \leq Q$. Again since R is right krull homogenous and N is right w.i.i ideal of R so using theorem (3.4) (ii) above we get that $N = P$. We prove now that R is also a left krull homogenous ring. Suppose this is not true. So assume $B \neq 0$. Then by proposition (1.3) of section (1) above we have that $l(B) \neq 0$, and $|R/l(B)|_l < |R|_l$. Note that by theorem (3.4) (i) above we have that $\Lambda' \leq \Lambda$ and $P \leq Q$. Since R/Q is a left krull homogenous semiprime ring with $|R/Q|_l = |R|_l$, so we must have that $l(B) \not\subset Q$. Since $P \leq Q$, so we must have that $l(B) \not\subset P$ as well. Since R/P is a right krull homogenous semiprime ring with $|R/P|_r = |R|_r$, so $(l(B) + P)/P \cap C(R/P) \neq \Phi$. It is not difficult to see that this implies that $l(B) \cap c(P) \neq \Phi$. Let $d \in l(B) \cap c(P)$. Note that by theorem (8) of [5] R has an artinian quotient ring. Hence, since $N = P$ and R has an artinian quotient ring, so using Small's theorem (see [6]) we get that $d \in l(B) \cap c(0)$. Therefore $dB = 0$, gives that $B = 0$, a contradiction to our original assumption that $B \neq 0$. Thus we must have that $B = 0$ and hence R is also a left krull homogenous ring.

**Theorem (3.6)** :- Let R be a noetherian ring with nilradical N. Let $\Lambda$, $\Lambda'$ be as in definition (1.4) above. Let $\Lambda \cap \Lambda' \neq \Phi$. If $V = \Lambda \cap \Lambda'$ let $Y = \{ \cap p \,/\, p \,\varepsilon\, V \}$ as in definition (3.3) above, then the following hold true;

(a) If N is right w.i.i ideal, then $V = \Lambda$, and $Y = P$ is a right localizable ideal of R.

(b) If N is left w.i.i ideal, then $V = \Lambda'$, and $Y = Q$ is a left localizable ideal of R.

**Proof :-** (a) Consider the factor ring R/A. Then R/A is a right krull homogenous ring with $|R/A|_r = |R|_r$. Since $V \neq \Phi$, thus by theorem (1.8) of section (1) above we must have that $|R/A|_l = |R|_l$. Since N is a right w.i.i ideal, so by theorem (3.5) proved above R/A is also a left krull homogenous ring. Hence by the main theorem of [7], namely, theorem (4) of [7], we must have that $\Lambda'(R/A) = \Lambda(R/A)$. Moreover, since by proposition (1.2) above, $\Lambda(R/A) = \Lambda(R)$, so we get that $\Lambda(R) = \Lambda'(R/A)$. By theorem (1.5) of section (1), since $|R/A|_l = |R|_l$, so we must have that $\Lambda'(R/A) \leq \Lambda'(R)$. Thus from the foregoing we finally get that $\Lambda(R) \leq \Lambda'(R)$. Hence we must have that $V = \Lambda(R)$ and hence $Y = P$, where $P = \{ \cap p \,/\, p \,\varepsilon\, \Lambda \}$, which is the right krull homogenous part of N. Thus by the main theorem of [9], namely theorem (2.10) of [9] (see theorem (3.2) of this paper) we get that P is right localizable. So Y is a right localizable ideal of R.

(b) The proof for (b) is the left analogue of (a). In this case $V = \Lambda'(R)$ and $Y = Q$, and hence $Y = Q$ is a left localizable ideal of R.

From the above theorem we get the main theorem of section (3) for the case $\Lambda \cap \Lambda' \neq \Phi$, which we state below;

**Theorem (3.7) :-** Let R be a noetherian ring with nilradical N. Let R be a non-disjoint ring. Then the following are equivalent conditions on R ;

(a) N is right and left w.i.i. ideal of R.

(b) P = Q is both a right and left localizable ideal of R.

**Proof :-** (a) => (b) :- Let $\Lambda, \Lambda'$ be as in definition (1.4) above. If R is a non-disjoint ring, that is if $\Lambda \cap \Lambda' \neq \Phi$ **then** set $V = \Lambda \cap \Lambda'$ and let $Y = \{ \cap p \, / \, p \, \varepsilon \, V \}$. Then since N is a right and a left w.i.i of R, so from theorem (3.6) above we have that $V = \Lambda = \Lambda'$, and Y = P = Q is a right and a left localizable ideal of R.

(b) => (a) follows from the main theorem of [9] ( see theorem (3.2) of this paper stated above ).

Combining theorem (2.9) and theorem (3.7) of this paper we get theorem (3.8) below which is also the main theorem of this paper ;

**Theorem (3.8) (Main theorem ) :-** Let R be a Noetherian ring with nilradical N. Let P and Q be the semiprime ideals that are the right and the left krull -homogenous parts of N. Then the following hold true ;

(a) If R is a disjoint ring, then the nilradical N of R is a right and a left weakly ideal invariant ideal of R. Hence N is a right and a left localizable semiprime ideal of R.

(b) If R is a non-disjoint ring then the following are equivalent conditions on R ;

(a) N is right and left w.i.i. ideal of R.

(b) P = Q is a right and left localizable ideal of R.

An interesting feature of a non-disjoint ring is the statement of the following corollary (3.9) given below ;

**Corollary (3.9) :-** Let R be a noetherian ring with nilradical N . Let $\Lambda , \Lambda'$ , A and B be as in definitions (1.4 ) and (1.1 ) above . Let $\Lambda \cap \Lambda' \neq \Phi$ . If $V = \Lambda \cap \Lambda'$ and let $Y = \{ \cap p / p \varepsilon V \}$ as in definition (3.3) above , then the following hold true ;

(a) N is a right and a left w.i.i. ideal of R implies that R is a weakly ideal krull symmetric ring ( that is A = B ) .

(b) A=B implies that $\Lambda = \Lambda'$ . Equivalently , this means that

a weakly ideal krull symmetric ring is a weakly krull symmetric ring .

**Proof :-** (a) Suppose first that N is a right w.i.i ideal of R . Then theorem (2.6) stated above guarantees that N(R/A) is a right w.i.i and thus theorem (3.5) stated above guarantees that R/A is a left krull –homogenous ring . From theorem (1.9) above we have that $| R/A |_l = | R |_l$ , hence we must have that $B \leq A$ . Similarly using the left w.i.i of N , we get that $A \leq B$ . Thus A = B .

(b) (b) is the same as the statement of theorem (2.9) of [8] .

**Remark :-** In connection with corollary (3.9) above we mention that converse questions of the statements in (a) and (b) of corollary (3.9) are open questions .

In general in the absence of the w.i.i of the nilradical of the ring R of the above theorem (3.7) we have the following result ;

**Theorem (3.10) :-** Let R be a noetherian ring with nilradical N . Let $\Lambda , \Lambda'$ be as in definition (1.4 ) above . If $\Lambda \cap \Lambda' \neq \Phi$ , then set $V = \Lambda \cap \Lambda'$ and $Y = \{ \cap p / p \varepsilon V \}$ , then the following hold true ;

(i) $A+B \leq Y$ and $V( R/(A+B) = \Lambda(R/A+B) = \Lambda'(R/A+B)$ .

(ii) Y is a right localizable semiprime ideal of R if and only if N(R/A+B), the nilradical of the ring R/A+B is a right w.i.i ideal of the ring R.

(iii) Y is a left localizable semiprime ideal of R if and only if N(R/A+B), the nilradical of the ring R/A+B is a left w.i.i ideal of the ring R.

**References:-**

(1) A.K. Boyle, E.H. Feller; Semi critical modules and K-primitive rings, "L.M.N"; No. 700, Springer Verlag

(2) K.A. Brown, T.H. Lenagan and J.T. Stafford, "Weak Ideal Invariance and localisation"; J. Lond. math.society (2), 21, 1980, 53-61.

(3) K.R. Goodearl and R.B. Warfield, "An Introduction to Non commutative Noetherain Rings", L.M.S., student texts, 16.

(4) R. Gordon, "some Aspects of Non Commutative Noetherian Rings", L.N.M.; vol.545; Spg. Vlg., 1975, 105-127

(5) G.Krause, T.H. Lenagan and J.T. Stafford, " Ideal Invariance and Artinian Quotient Rings "; Journal of Algebra 55, 145-154 (1978).

(6) L.W.Small; Orders In Artinian Rings, J.Algebra 4 (1966). 13-41. MR 34, #199.

(7) C.L. Wangneo ; On A Certain Krull Symmetry Of a Noetherian Ring. arXiv:1110.3175, v2 [math.RA], Dec.2011.

(8) C.L. Wangneo ; On the weak Krull Symmetry Of a Noetherian Ring. arXiv: 1511.02678v2 [math.RA], April 2,2016.

(9) C.L. Wangneo ; On some conditions on a noetherian ring arXiv: :1503.00462v5 [math.RA], April 2, 2016